\documentclass[12]{amsart}
\usepackage{amssymb,mathdots}
\usepackage{mathrsfs}
\usepackage[applemac]{inputenc}
\usepackage{amsmath, graphicx}
\usepackage{wrapfig}
\usepackage{enumerate}
\usepackage{url}

\def\scaledpicture#1by#2(#3scaled#4){{
\dimen0=#1  \dimen1=#2
\divide\dimen0 by 1000 \multiply\dimen0 by #4
\divide\dimen1 by 1000 \multiply\dimen1 by #4
\picture \dimen0 by \dimen1 (#3 scaled #4)}}
\def\dfigure#1by#2(#3scaled#4offset#5:#6)
    {\medskip
     \vglue 2mm minus 2mm
     $$
       \hbox{
         \hglue#5
         {\scaledpicture #1 by #2 (#3 scaled #4)}
       }
     $$
     \par\goodbreak
     \vglue 2mm minus 2mm
     \medskip}

\def\qmod#1#2{{\hbox{}^{\displaystyle{#1}}}\!\big/\!\hbox{}_{
\displaystyle{#2}}}

\def\resto#1#2{{
#1\hskip 0.4ex\vline_{\hskip 0.2ex\raisebox{-0,2ex}
{{${\scriptstyle #2}$}}}}}

\def\C{{\mathbb C}}

\def\N{{\mathbb N}}

\def\P{{\mathbb P}}

\def\R{{\mathbb R}}
\def\Z{{\mathbb Z}}

\def\map{\longrightarrow}
\def\textmap#1{\mathop{\vbox{\ialign{
                                  ##\crcr
      ${\scriptstyle\hfil\;\;#1\;\;\hfil}$\crcr
      \noalign{\kern 1pt\nointerlineskip}
      \rightarrowfill\crcr}}\;}}
\def\bigtextmap#1{\mathop{\vbox{\ialign{
                                  ##\crcr
      ${\hfil\;\;#1\;\;\hfil}$\crcr
      \noalign{\kern 1pt\nointerlineskip}
      \rightarrowfill\crcr}}\;}}
      
\newcommand{\cal}{\mathcal}
\def\textlmap#1{\mathop{\vbox{\ialign{
                                  ##\crcr
      ${\scriptstyle\hfil\;\;#1\;\;\hfil}$\crcr
      \noalign{\kern-1pt\nointerlineskip}
      \leftarrowfill\crcr}}\;}}

\def\g{{\mathfrak g}}

\newtheorem{sz}{Satz}
\newtheorem{thry}[sz]{Theorem}
\newtheorem{pr}[sz]{Proposition}

\newtheorem{dt}[sz]{Definition}

\begin{document}
 
\def\tr{\mathrm {Tr}}
\def\End{\mathrm {End}}
\def\Aut{\mathrm {Aut}}
\def\Spin{\mathrm {Spin}}
\def\U{\mathrm{U}}
\def\SU{\mathrm {SU}}
\def\SO{\mathrm {SO}}
\def\PU{\mathrm {PU}}
\def\GL{\mathrm {GL}}
\def\spin{\mathrm {spin}}
\def\u{\mathrm {u}}
\def\su{\mathrm {su}}
\def\so{\mathrm {so}}
\def\pu{\mathrm {pu}}
\def\Pic{\mathrm {Pic}}
\def\Iso{\mathrm {Iso}}
\def\NS{\mathrm{NS}}
\def\deg{\mathrm {deg}}
\def\Hom{\mathrm{Hom}}
\def\Herm{\mathrm{Herm}}
\def\Vol{{\rm Vol}}
\def\pf{{\bf Proof: }}
\def\id{ \mathrm{id}}
\def\Im{\mathrm{Im}}
\def\im{\mathrm{im}}
\def\rk{\mathrm {rk}}
\def\ad{\mathrm {ad}}
\def\spc{\mathrm{Spin}^c}
\def\U2{\mathrm{U(2)}}
\def\niq{=\kern-.18cm /\kern.08cm}
\def\Ad{\mathrm {Ad}}
\def\RSU{\R\mathrm{SU}}
\def\ad{{\rm ad}}
\def\dva{\bar\partial_A}
\def\da{\partial_A}
\def\p{{\rm p}}
\def\sp{\Sigma^{+}}
\def\sm{\Sigma^{-}}
\def\spm{\Sigma^{\pm}}
\def\smp{\Sigma^{\mp}}
\def\oo{{\scriptstyle{\cal O}}}
\def\ooo{{\scriptscriptstyle{\cal O}}}
\def\sw{Seiberg-Witten }
\def\pa{\partial_A\bar\partial_A}
\def\Dr{{\raisebox{0.17ex}{$\not$}}{\hskip -1pt {D}}}
\def\gr{{\scriptscriptstyle|}\hskip -4pt{\g}}
\def\subsetint{{\  {\subset}\hskip -2.45mm{\raisebox{.28ex}
{$\scriptscriptstyle\subset$}}\ }}

\title
{Real determinant line bundles}
\author{Christian Okonek \& Andrei Teleman}
 \begin{abstract}   This article is  an expanded  version of the talk given by \hbox{Ch.\hspace{0.7mm}O.} at the  Second Latin Congress on "Symmetries in
Geometry and Physics" in Curitiba, Brazil in December 2010. In this version we explain the topological and gauge-theoretical aspects of our paper "{\it Abelian Yang-Mills theory on Real tori and Theta divisors of Klein surfaces}" \cite{OT3}.
 \end{abstract} 
\maketitle

\section{Introduction}
 
In an ongoing research project  we intend to develop a version of Seiberg-Witten theory and gauge theoretic Gromov-Witten theory in the presence of Real structures. The main problem is that the  corresponding virtual fundamental class will be a homology  class with coefficients in a  local  system (i.e., a locally constant sheaf) which, in order to construct  $\Z$-valued  invariants, must  be determined explicitly.
\vspace{2mm}\\ 
{\bf Example:}
Let $(C,\iota)$ be a Klein surface, i.e., a Riemann surface $C$ endowed with an anti-holomorphic involution $\iota$. Put $g:=h^0(\omega_C)$  and suppose   $C^\iota\ne\emptyset$. 

\vspace{2mm}
\begin{center} 
\includegraphics[width=4cm]{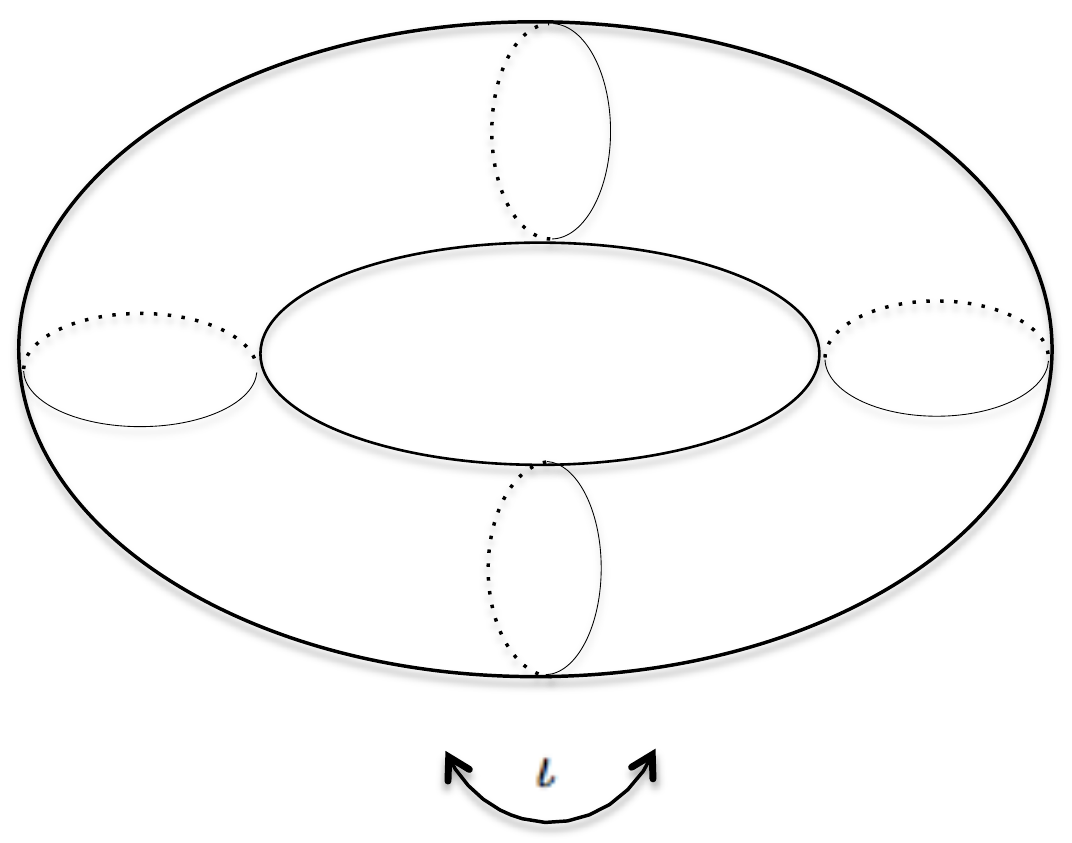}
\end{center}
For every $d\in\N$ we obtain   induced Real structures $\hat\iota:\Pic^d(C)\to\Pic^d(C)$, and $\tilde\iota: C^{(d)}\to C^{(d)}$ on the complex manifolds $\Pic^d(C)$, and $C^{(d)}$. 

In the example illustrated   in the picture above, we have $(C^{(2)})^{\tilde\iota}\simeq T^2{\textstyle\coprod }(\P^2_\R\#\P^2_\R)$. 
The  fundamental class of the manifold $(C^{(d)})^{\tilde\iota}$   is an element 
$$[(C^{(d)})^{\tilde\iota}]\in H_d((C^{(d)})^{\tilde\iota},\oo_{(C^{(d)})^{\tilde\iota}} )\ ,$$
where $\oo_{(C^{(d)})^{\tilde\iota}}$ denotes its orientation sheaf. To understand this coefficient system we use the cartesian diagram:
$$
\begin{array}{c}
\unitlength=1mm
\begin{picture}(32,12)(-4,-4)
\put(-10,4){$C\times C^{(d)}$}
\put(6,5){\vector(2,0){13}}
\put(22,4){$ C^{(d)}$}
\put(-3,1){\vector(0, -3){8}}
\put(-12, -12){$C\times\Pic^d(C)$}
\put(-2,-3){$\id\times\rho$}
\put(21,-3){$\rho$}
\put(13,7){$\tilde q$}
\put(13,-9){$ q$}
\put(24,1){\vector(0, -3){8}}
\put(20, -12){$\Pic^d(C)$}
\put(10,-11){\vector(2,0){8}}
\end{picture} 
\end{array}\  
$$
\vspace{6mm}\\
One has the following effective divisors:
\begin{enumerate}
\item  {\it the universal divisor} $\Delta:=\{(x,\delta)\in C\times C^{(d)}|\ x\in\delta\}$ on the product $C\times C^{(d)}$;
\item the ample divisor  
$D_{x_0}:=\{\delta\in C^{(d)}|\ x_0\in \delta\}$ on $C^{(d)}$, where  $x_0\in C$ is a fixed point.
\end{enumerate}
We denote by ${\cal O}(\Delta)$, respectively ${\cal O}(D_{x_0})$ the corresponding line bundles on the manifolds $C\times C^{(d)}$,  $C^{(d)}$ respectively. On the other hand, on $C\times\Pic^d(C)$ we consider the Poincaré line bundle with respect to $x_0$, denoted ${\cal P}_{x_0}$, which is characterized up to isomorphism by the properties:
$$\resto{{\cal P}_{x_0}}{C\times\{[{\cal L}]\}}\simeq {\cal L}\ \ \forall[{\cal L}]\in\Pic^d(C)\ ,\ \resto{{\cal P}_{x_0}}{\{x_0\}\times\Pic^d(C)}\simeq{\cal O}_{\Pic^d(C)}\ .
$$
We have a universal short exact sequence on $C\times C^{(d)}$:
$$0\map {\cal O}\map {\cal O}(\Delta)\map {\cal O}_\Delta(\Delta)\map 0\  
$$
Using the well known facts \cite{ACGH}
\begin{enumerate}
\item[i)] $T_{C^{(d)}}= \tilde q_*{\cal O}_\Delta(\Delta)$,
\item[ii)] ${\cal O}(\Delta)=(1\times\rho)^*{\cal P}_{x_0}\otimes \tilde q^*{\cal O}(D_{x_0})$,
\end{enumerate}
and taking into account that determinant line bundles  commute with base change \cite{KM}, we obtain:
\begin{equation}\label{T}\det T_{C^{(d)}}=\det(\tilde q_*{\cal O}_\Delta(\Delta))=\det(\tilde q_!{\cal O}_\Delta(\Delta))=\det(\tilde q_!{\cal O}(\Delta))\otimes \det(\tilde q_!{\cal O})^\vee\simeq
$$
$$\simeq \det(\tilde q_![(1\times\rho)^*{\cal P}_{x_0}\otimes\tilde q^*{\cal O}(D_{x_0})])\simeq \rho^*[\det(q_!{\cal P}_{x_0})]\otimes{\cal O}((d+1-g)D_{x_0})\  
\end{equation}

When $x_0\in C^{\iota}$ all  bundles in these formulae are naturally Real holomorphic line bundles, and the isomorphisms in   (\ref{T}) are isomorphisms of Real line bundles. In particular we see that the orientation  line bundle $\det T_{(C^{(d)})^{\tilde\iota}}$ of the manifold $(C^{(d)})^{\tilde\iota}$ -- which  can be identified with the fixed point locus of $\det T_{C^{(d)}}$  --  is determined by the numerical data $g$, $d$ and the Real holomorphic line bundle $\det(q_!{\cal P}_{x_0})$ on $(\Pic^d(C),\hat\iota)$. Therefore, in order to determine the orientation sheaf $\oo_{(C^{(d)})^{\tilde\iota}}$, we have to understand this Real determinant line bundle. To this end consider the Real isomorphism $\varphi_{x_0}:\Pic^0(C)\map\Pic^d(C)$ given by
$$\varphi_{x_0}([{\cal L}])=[{\cal L}\otimes{\cal O}_C(dx_0)]\ .
$$
\begin{pr} Let $\Theta:=\{[{\cal L}]\in\Pic^{g-1}(C)|\ h^0({\cal L})>0\}$ be the geometric theta divisor on $\Pic^{g-1}(C)$. There exists a canonical isomorphism of Real holomorphic line bundles
$$\varphi_{x_0}^*[\det(q_!{\cal P}_{x_0})]\simeq{\cal O}_{\Pic^0(C)}(\Theta-[{\cal O}_C((g-1)x_0)])\ .$$
\end{pr}
{\ }\\
{\bf Remark:} \begin{enumerate}
\item[i)] The isomorphism type of the underlying differentiable {\it complex} line bundle of the right hand side can be computed using the Grothendieck-Riemann-Roch theorem for proper morphisms, and the result is independent of $x_0$. On the other  hand, we will see that the isomorphism type of the  {\it Real} structure of this line bundle does depend on $x_0$, more precisely on the connected component of $x_0$ in $C^\iota$. This striking result shows that there cannot exist any Grothendieck-Riemann-Roch-type  theorem to compute this Real isomorphism type.
\item[ii)] When one considers more generally quot-schemes $\mathrm{Quot}^d_{\cal E}$ instead of the symmetric power $C^{(d)}=\mathrm{Quot}^d_{\cal O}$, one has to deal with virtual fundamental classes in order to define in a coherent way Real gauge theoretical Gromov-Witten invariants \cite{OT2}.
\end{enumerate}

Therefore, our aim is now to determine ${\cal O}_{\Pic^0(C)}(\Theta-[{\cal O}_C((g-1)x_0)])$ as a Real holomorphic line bundle on the Real torus $(\Pic^0(C),\hat\iota)$.
\section{Real line bundles}

Let $(X,\tau)$ be a Real  space, i.e. a CW space endowed with an involution. A Real line bundle (in the sense of Atiyah \cite{A1}) over $(X,\tau)$ is a pair $(L,\tilde\tau)$, where $L\to X$ is a complex line bundle, and $\tilde\tau:L\to L$ is a fibrewise anti-linear  isomorphism over $\tau$ such that $\tilde\tau^2=\id_L$. The isomorphism classes of Real line bundles over  $(X,\tau)$ form a group, which is naturally isomorphic to the Grothendieck cohomology group \cite{G}
$$H^1_{\Z_2}(X,\underline{S}^1(1))= H^2_{\Z_2}(X,\underline{\Z}(1))\simeq[(X,\tau),(\P^\infty(\C),\bar{\ })]^{\Z_2}\ ,
$$
where $\underline{S}^1(1)$ denotes the $\Z_2$-sheaf of germs of continuous $S^1$-valued functions on $X$ endowed with the  $\Z_2$-action defined by the involution induced by $\tau$ and conjugation on $S^1$, and $\underline{\Z}(1)$ denotes the  $\Z_2$-sheaf  with fibre $\Z$ endowed with the  $\Z_2$-action defined by the involution induced by $\tau$ and $-\id_\Z$ (see  \cite{G} for the cohomology theory of equivariant sheaves) .

\begin{pr} \label{CanExSeq} Suppose that $X^\tau\ne\emptyset$. There exists a canonical exact sequence
$$0\to H^1_{\Z_2}(H^1(X,\Z)(1))\to H^2_{\Z_2}(X,\underline{\Z}(1))\textmap{c_1} H^2(X,\Z)^{-\tau^*}\textmap{\ooo} H^2_{\Z_2}(H^1(X,\Z)(1))\ .
$$
\end{pr}
We also have a natural restriction map
$$H^1_{\Z_2}(X,\underline{S}^1(1))=H^2_{\Z_2}(X,\underline{\Z}(1))\map H^2(X^\tau,
\underline{\Z}(1))=H^1(X^\tau,\Z_2)\ ,
$$
which maps the isomorphism class $[L,\tilde\tau]$ of a Real line bundle  on $(X,\tau)$ onto $w_1(L^{\tilde\tau})$. Here $L^{\tilde\tau}$ is regarded as a real line bundle on the fixed point locus $X^\tau$.
\vspace{3mm}\\
{\bf Example:} For a Klein surface $(C,\iota)$ one has
$$H^2_{\Z_2}(C,\underline{\Z}(1))\simeq H^2(C,\Z)\times_{\Z_2} H^1(C^\iota,\Z_2) \ ,
$$
the isomorphism being given by $[L,\tilde\iota]\mapsto (c_1(L),w_1(L^{\tilde\iota}))$. The fact that the pair on the right belongs to the fibre product $H^2(C,\Z)\times_{\Z_2} H^1(C^\iota,\Z_2)$ follows from the identity
$$\langle c_1(L),[C]\rangle\equiv \langle w_1(L^{\tilde\iota}),[C^\iota]\rangle\ (\hbox{mod } 2)\ .
$$
\section{Gauge theoretical relevance}

Let $(X,\tau)$ be a compact Real Riemannian manifold (i.e., $\tau$ is an isometric involution of $X$) 	and let $L$ be a Hermitian line bundle on $X$. We denote by ${\cal T}(L)$ the moduli space of Yang-Mills connections, and  put
$${\cal T}_X:=\coprod_{c\in H^2(X,\Z)} {\cal T}(L_c)\ ,$$
where $L_c$ is a Hermitian line bundle on $X$  with $c_1(L_c)=c$. Note that ${\cal T}(L_c)$ depends only on $c$ (i.e., is independent of the choice of $L_c$) up to canonical isomorphism. The manifold ${\cal T}_X$ comes with a natural  Real structure  $\hat\tau:{\cal T}_X\to {\cal T}_X$ defined by pull-back of Yang-Mills connections.
\begin{pr} \label{Ff} Suppose $X^\tau\ne\emptyset$. 
\begin{enumerate}
\item The following conditions are equivalent:
\begin{enumerate}
\item ${\cal T}(L_c)^{\hat\tau}\ne\emptyset$,
\item $\tau^*(c)=-c$, $\ooo(c)=0$,
\item $L_c$ admits Real structures with respect to the fixed involution $\tau$ on $X$.
\end{enumerate}
\item There exists a  natural morphism
$$F:{\cal T}_X^{\hat\tau}\map H^2_{\Z_2}(X,\underline{\Z}(1))\ ,
$$
which induces an isomorphism $f:\pi_0({\cal T}_X^{\hat\tau})\map H^2_{\Z_2}(X,\underline{\Z}(1))$.
\end{enumerate}
\end{pr}
{\ }\\
{\bf Problem:} Compute  the group $H^2_{\Z_2}(X,\underline{\Z}(1))$ classifying Real line bundles on a given Real 4-manifold $(X,\tau)$. Note that this is taken care of by Proposition \ref{CanExSeq} when the first Betti number of $X$ vanishes.
\section{Abelian Yang-Mills connections on a torus}
Let $V$ be an Euclidian vector space, $\Lambda\subset V$ a maximal lattice, and $T:=V/\Lambda$ the associated flat torus. Recall that we have   canonical isomorphisms $H^k(T,\Z)=\mathrm{Alt}^k(\Lambda,\Z)$. Fix $u\in H^2(T,\Z)$, denote by the same symbol the corresponding anti-symmetric bilinear form $u:\Lambda\times\Lambda\to\Z$, and by $u_\R$ its $\R$-linear extension, which can  also be  regarded as a harmonic 2-form on $T$. Using the notations of the previous section we have
$${\cal T}(L_u)=\left\{[A]\in\qmod{{\cal A}(L_u)}{\cal G}\ \vline\ \frac{i}{2\pi} F_A=u_\R\right\}\ .
$$
\begin{dt} A $u$-character on $\Lambda$ is a map $\alpha:\Lambda\to S^1$ such that
$$\alpha(\lambda+\lambda')=\alpha(\lambda)\alpha(\lambda')e^{\pi i u(\lambda,\lambda')}\  \ \forall\lambda,  \lambda'\in\Lambda\ .
$$
\end{dt}
Note that the set $\Hom_u(\Lambda,S^1)$ of $u$-characters on $\Lambda$ has a natural $\Hom(\Lambda,S^1)$-torsor structure.

Associated to any element $\lambda\in\Lambda$ one has a {\it standard loop} $c_\lambda:[0,1]\to  V/\Lambda$ given by $c_\lambda(t):=[\lambda t]$.
\begin{pr} Taking  holonomy along  standard loops induces a bijection
$$h:{\cal T}(L_u)\to \Hom_u(\Lambda,S^1)\ ,
$$
which maps a Yang-Mills  class $[A]\in {\cal T}(L_u)$ to the $u$-character $\lambda\mapsto \bar h^A_{c_\lambda}$.
\end{pr}
{\ }\vspace{-3mm}\\
{\it Idea:} The  holonomy along the boundary of a 2-simplex is determined by the curvature form. We apply this principle  to the image in $T$ of a 2-simplex $[0,\lambda,\lambda']\subset V$,  and see that the holonomy of a Yang-Mills connection along the standard loops defines a $u$-character. Letting $u$ vary in $\mathrm{Alt}^2(\Lambda,\Z)$ we obtain
\begin{thry} \label{AH} (Appel-Humbert theorem for Abelian Yang-Mills connections). The holonomy map defines a canonical isomorphism
$$h:{\cal T}_{T}\map \coprod_{u\in\mathrm{Alt}^2(\Lambda,\Z)}\Hom_u(\Lambda,S^1)\ .
$$
\end{thry}
{\ }\\
{\bf Remark:} The classical Appel-Humbert theorem, describing holomorphic line bundles on complex tori, follows from this result and the Kobayashi-Hitchin correspondence \cite{LT}: if $J$ is a compatible complex structure on $V$, we define a Hermitian form $H_u$ on  the complex vector space $(V,J)$  by
$$H_u(v,w):= u_\R(v,Jw)+i u_\R(v,w)\ .
$$
Clearly one has $\im H_u=u_\R$. We obtain a commutative diagram
\vspace{1mm}\\
$$
\begin{array}{c}
\unitlength=1mm
\begin{picture}(32,12)(-8,-4)
\put(-36,4){$\displaystyle\coprod_{\begin{array}{c} 
\scriptstyle H\ \mathrm{Hermitian}\vspace{-1mm}\\ \scriptstyle\im H(\Lambda\times\Lambda)\subset\Z \end{array}}
\hspace{-7mm} \Hom_{\im H}(\Lambda,S^1)$}
\put(5,5.5){\line(2,0){10}}
\put(5,4.5){\line(2,0){10}}
\put(12,4){$\displaystyle\coprod_{\begin{array}{c} 
\scriptstyle u\in\mathrm{Alt}^2(\Lambda,\Z)\vspace{-1mm}\\ \scriptstyle  J^*(u_\R)=u_\R \end{array}}\hspace{-7mm} \Hom_{u}(\Lambda,S^1)$}
\put(-22,-7){\vector(0, -3){12}}
\put(-28, -24){$\Pic(T,J)$}
\put(-20,-13){$AH$}
\put(16.5,-13){$h^{-1}$}
\put(-2,-21){$KH$}
\put(23,-7){\vector(0, -3){12}}
\put(16, -24){$\displaystyle\coprod_{u\in \mathrm{NS}(T,J)}\hspace{-5mm}{\cal T}(L_u)$}
\put(16,-23){\vector(-2,0){26}}
\put(42,-23){,}
\end{picture} 
\end{array} 
$$
\vspace{26mm}\\
where the map $KH$ is the Kobayashi-Hitchin correspondence between equivalence classes of HE connections and isomorphism classes of polystable holomorphic bundles. We recall that a Hermitian connection on a differentiable Hermitian  bundle over a compact Kähler manifold is Hermitian-Einstein if and only if it is Yang-Mills and its curvature has type (1,1).
\vspace{4mm}\\
{\bf Remark:} Theorem \ref{AH} gives an interesting geometric interpretation of the classical Appel-Humbert data $(H,\alpha)$ as curvature, respectively holonomy of Hermitian-Einstein connections.
\vspace{4mm}\\
{\bf Problem 2:} Generalize Theorem \ref{AH} to higher harmonic Deligne cohomology, i.e., higher Abelian gerbes with harmonic curvature. For the definition of these spaces we refer to \cite{DK}. 
\section{Real line bundles on Real tori}\label{RLBT}

Let $V$ be an Euclidian vector space, $\tau:V\to V$ a linear isometric involution, and $\Lambda\subset V$ a $\tau$-invariant maximal lattice. We denote by the same symbol $\tau$ the induced involution on the torus $T:=V/\Lambda$, and by $\hat\tau:{\cal T}_T\to{\cal T}_T$ the induced involution on the total moduli space ${\cal T}_T$ of Yang-Mills connections on Hermitian line bundles over $T$. Using the notations and the results of the previous section we see that  $\hat\tau$ acts -- via the AH description given by Theorem \ref{AH} --  by the formula
$$\hat \tau ( u,\alpha)=(-\tau^* u,\overline{\tau^*\alpha})\ ,
$$
which shows that a fixed point $(u,\alpha)$ of $\hat\tau$ must satisfy the equations
\begin{enumerate}
\item $\resto{\alpha}{\Lambda^\tau}\in\Hom(\Lambda^\tau,\{\pm 1\})$,
\item $\alpha(\lambda+\tau\lambda)=e^{\pi i u(\lambda,\tau\lambda)}$ $\forall\lambda\in\Lambda$.
\end{enumerate} 
Therefore, restricting the map $h$ of Theorem \ref{AH} to the fixed point locus of $\hat\tau$ we get a map
$$h^\tau:{\cal T}_T^{\hat\tau}\map \mathrm{Alt}^2(\Lambda,\Z)^{-\tau^*}\times_{\Hom((\id+\tau)\Lambda,\Z_2)}\Hom(\Lambda^\tau,\Z_2)
$$
given by $h^\tau(u,\alpha):=(u,\resto{\alpha}{\Lambda^\tau})$. Note that Proposition \ref{Ff} yields a natural bijection $f:\pi_0({\cal T}_T^{\hat\tau})\to H^2_{\Z_2}(T,\underline{\Z}(1))$; clearly $h^\tau$ factorizes through $f$. Note also that one has canonical identifications
$$\mathrm{Alt}^2(\Lambda,\Z)^{-\tau^*}=H^2(T,\Z)^{-\tau^*}\ ,\ \Hom(\Lambda^\tau,\Z_2)=H^1(T^\tau_0,\Z_2)\ ,
$$
where $T^\tau_0=V^\tau/\Lambda^\tau$ is the connected component of 0 in the fixed point locus $T^\tau$.
\begin{thry} Let $(T,\tau)$ be a Real torus with $T^\tau\ne\emptyset$. The mixed characteristic class $(c_1,w_1)$ induces an isomorphism
$$cw:H^2_{\Z_2}(T,\underline{\Z}(1))\map H^2(T,\Z)^{-\tau^*}\times_{H^1(T,\Z_2)} H^1(T^\tau_0,\Z_2)
$$
given by $cw([L,\tilde\tau])= (c_1(L),w_1(\resto{L^{\tilde\tau}}{T^\tau_0}))$.
\end{thry}
{\ }\\
{\bf Remark:} The essential facts used in the proof are:
\begin{enumerate}
\item The second component $\alpha$ of an Appel-Humbert datum describes the holonomy of the corresponding Yang-Mills connection,
\item The connected component decomposition of the fixed point locus $T^\tau$ has the form
$$T^\tau=\coprod_{[\mu]\in \frac{1}{2}\Lambda^{-\tau}/\frac{1}{2}(\id-\tau)\Lambda} T^\tau_{[\mu]}\ ,
$$
where $T^\tau_{[\mu]}:=T_0^\tau+[\mu]$; the Stiefel-Whitney class $w_1(L^{\tilde\tau})\in H^1(T^\tau,\Z_2)$ is determined by the Stiefel-Whitney class $w_1(\resto{L^{\tilde\tau}}{T^\tau_0})$ of the restriction of $ {L^{\tilde\tau}}$ to $T^\tau_0$. This follows from the difference formula
$$w_1(\resto{L^{\tilde\tau}}{T^\tau_{[\mu]}})-w_1(\resto{L^{\tilde\tau}}{T^\tau_0})=  {u(2\mu,\cdot)}\ (\mathrm{mod}\ 2)\ .
$$
\end{enumerate}
{\bf Problem 3:} Silhol \cite{S2} constructed moduli spaces of Real Abelian varieties $(A,\tau)$ endowed with a compatible principal polarization $c\in H^2(A,\Z)^{-\tau^*}$. For certain questions it is more natural to have moduli spaces of Real principally polarized Real Abelian varieties, i.e., triples $(A,\tau,cw)$, where $cw\in H^2_{\Z_2}(A,\underline{\Z}(1))$ is a class defining a principal polarization. These   moduli spaces will be  finite covers of Silhol's spaces.
\section{Determinant bundles of Klein surfaces}\label{DBKS}
Let $(C,\iota)$ be a Klein surface,   $C^\iota=\coprod_{i=1}^nC_i$  the connected component decomposition of the fixed point locus $C^\iota$, and  let  $x_0\in C^\iota$.  
Our goal is to determine  the topological type of the Real line bundle 
$${\cal L}_{x_0}:={\cal O}_{\Pic^0(C)}(\Theta-[{\cal O}_C((g-1)x_0)])\ .$$
Taking into account the results explained in the previous section, this topological type is determined by the mixed characteristic class
$$
cw({\cal L}_{x_0})\in \mathrm{Alt}^2(H^1(X,\Z),\Z)^{({\iota^*})^*}\times_{\Hom((\id-\iota^*)H^1(C,\Z),\Z_2 )} \Hom(H^1(C,\Z)^{-\iota^*},\Z_2)\ .$$
The Grothendieck-Riemann-Roch theorem  allows us to identify the first component of the pair $cw({\cal L}_{x_0})$ \cite{ACGH};  the result is $c_1({\cal L}_{x_0})=u_C$, where  
$$u_C\in H^2(\Pic^0(C),\Z)=\mathrm{Alt}^2(H^1(C,\Z),\Z) $$
 is  the cup form:
$$u_C:H^1(C,\Z)\times H^1(C,\Z)\to\Z\ ,\   (\lambda,\lambda')\mapsto \langle\lambda\cup\lambda',[C]\rangle $$
Therefore, it suffices to determine explicitly  the Stiefel-Whitney class
\begin{equation}\label{w1}w_1(\resto{{\cal L}_{x_0}^{\tilde\iota}}{\Pic^0(C)^{\hat\iota}_0}): H^1(C,\Z)^{-\iota^*}\map\Z\ .
\end{equation}

Using topological arguments \cite{CN}, one can show that $H^1(C,\Z)^{-\iota^*}$ is generated by the subgroup $(\id-\iota^*) H^1(C,\Z)$ and the classes $[C_1]^\vee,\dots,[C_n]^\vee$, where $[C_i]^\vee$ denotes the Poincaré dual of the 1-homology class $[C_i]$ defined by an arbitrary orientation of the circle $C_i$. Since the restriction  of $w_1(\resto{{\cal L}_{x_0}^{\tilde\iota}}{\Pic^0(C)^{\hat\iota}_0})$ to $(\id-\iota^*) H^1(C,\Z)$ is known (it is determined by the Chern class $u_C$, see section \ref{RLBT}) we conclude that it suffices to compute the values of (\ref{w1}) on the classes $[C_i]^\vee$.\vspace{2mm}\\
We will proceed in two steps:
\begin{enumerate}
\item compute $cw({\cal L}_{[\kappa]})$, where ${\cal L}_{[\kappa]}$ is  {\it a symmetric theta line bundle},
\item compare $cw({\cal L}_{x_0})$ to $cw({\cal L}_{[\kappa]})$.
\end{enumerate}
\vspace{2mm}
(1) We recall that a theta characteristic is square root of $[\omega_C]$. Denote by $\theta(C)$ the set of theta characteristics of $C$, i.e.,
$$\theta(C):=\{[\kappa]\in\Pic^{g-1}(C)|\ \kappa^{\otimes 2}\simeq\omega_C\}\ .
$$
The cardinality of this set is $2^{2g}$, and for every   $[\kappa]\in\theta(C)$  we have an associated Mumford theta form $q_{[\kappa]}:\Pic^0(C)_2\map\Z_2$, defined  on the 2-torsion subgroup $\Pic^0(C)_2\subset\Pic^0(C)$ by
$$q_{[\kappa]}([\eta]):=h^0(\eta\otimes\kappa)-h^0(\kappa)\hbox{ (mod }2) 
$$
(see \cite{ACGH}). Using the natural identification 
$$\Pic^0(C)_2= \qmod{\frac{1}{2}H^1(C,\Z)}{H^1(C,\Z)}=H_1(C,\Z_2)\ ,
$$
we obtain a form $q_{[\kappa]}:H_1(C,\Z_2)\to\Z_2$ satisfying the {\it Riemann-Mumford} relations:
\begin{equation}\label{RM} q_{[\kappa]}(a+b)=q_{a}+q_{a}+ a\cdot b
\end{equation}
The main idea is to use the translation by a theta characteristic  instead of translation by $[{\cal O}_C((g-1)x_0)]$ to identify $\Pic^{g-1}$ with $\Pic^0(C)$. More precisely, we define 
$${\cal L}_{[\kappa]}:={\cal O}_{\Pic^0(C)}(\Theta-[\kappa])\ .
$$
\begin{thry} \label{AHLK} Let $[\kappa]\in \theta(C)$. Then the Appel-Humbert data of the holomorphic line bundle ${\cal L}_{[\kappa]}$ is $(u_C,\alpha_{[\kappa]})$, where $\alpha_{[\kappa]}:H^1(C,\Z)\to S^1$ is defined by
$$\alpha_{[\kappa]}(\lambda):=(-1)^{q_{[\kappa]}(\overline{\lambda\cap[C]})}\ .
$$
Here $\lambda\cap[C]\in H_1(C,\Z)$ is the Poincaré  dual of $\lambda$ and $\overline{\lambda\cap[C]})$ denotes  its image in $H_1(C,\Z_2)$.
\end{thry}
{\ }\\
{\it Idea of proof:} Since ${\cal L}_{[\kappa]}$ is   {\it symmetric} in the sense that $(-1)^*({\cal L}_{[\kappa]})\simeq {\cal L}_{[\kappa]}$, it follows that
$$\alpha_{[\kappa]}(\lambda)=(-1)^{\mathrm{mult}_{[\frac{1}{2}\lambda]}(\Theta-[\kappa])- \mathrm{mult}_{[0]}(\Theta-[\kappa])}
$$
(see \cite{BL}). Now we use {\it Riemann's singularity theorem}, which states
$$\mathrm{mult}_{[{\cal L}]}\Theta=h^0({\cal L})\ .
$$
\vspace{3mm}\\
Note that the $u_C$-character $\alpha_{[\kappa]}$ given by Theorem  \ref{AHLK} involves the algebraic geometric data $q_{[\kappa]}$. Therefore, we have to make an additional step, which will give a purely topological interpretation of the values $q_{[\kappa]}([C_i]_2)$ when $\hat\iota[\kappa]=[\kappa]$.
\begin{pr} Suppose $[\kappa]\in \theta(C)^{\hat\iota}$. Then
$$q_{[\kappa]}([C_i]_2)=\langle w_1(\kappa^{\tilde\iota}),[C_i]_2\rangle+1\ .
$$
\end{pr}
{\ }\\
{\it Idea of proof:}
We use the diagram
$$
\begin{array}{c}
\unitlength=1mm
\begin{picture}(32,12)(-4,-4)
\put(-7,4){$\theta(C)$}
\put(2,3.5){\vector(2,-1){6}}

\put(-3,1){\vector(0, -3){8}}
\put(-8, -12){$\mathrm{Spin}(C)$}
\put(-2,-3){$\xi$}
\put(5,3){$ q$}
\put(5,-9){$\omega$}
\put(8, -3.5){$Q(H_1(C,\Z_2),\cdot)$}
\put(2,-8){\vector(2,1){6}}
\put(33,-4){,}
\end{picture} 
\end{array}\  
$$
\vspace{8mm}\\
where $\mathrm{Spin}(C)$ denotes the set of isomorphism classes of $\mathrm{Spin}$-structures on $C$, $Q(H_1(C,\Z_2),\cdot)$ is the set of maps $H_1(C,\Z_2)\to\Z_2$ satisfying the Riemann-Mumford relations (\ref{RM}), $q$ is the assignment $[\kappa]\mapsto q_{[\kappa]}$  given by the Mumford theta form, $\xi$ is the bijection defined by Atiyah \cite{A2}, and $\omega$ is a bijection defined  by Johnson \cite{J} in purely topological terms.  We know by Mumford that $q$ is a morphism of $H^1(C,\Z_2)$-torsors, according to Atiyah \cite{A2} $\xi$ is a bijection,  by Johnson \cite{J} $\omega$ is a bijection,  and according to Libgober \cite{L} the diagram commutes. Combining these results, and using a direct geometrical argument, we obtain the following formula: 
\begin{equation}\label{OT}
\omega_{\xi_\kappa}([C_i]_2)=\langle w_1(\kappa^{\tilde\iota}),[C_i]_2\rangle +1 \  
\end{equation}
This completes step (1); for step (2)  we use the fact that every component of $\Pic^{g-1}(C)^{\hat\iota}$ contains $2^g$   Real theta characteristics. Concluding, we get the following  explicit formula  which, combined with the results of section \ref{DBKS}, completes the computation of the topological type of the Real line bundles ${\cal L}_{x_0}$.
\begin{pr}  
$$w_1(\resto{{\cal L}_{x_0}^{\tilde\iota}}{\Pic^0(C)_0})[C_i]_2^\vee=\left\{
\begin{array}{ccc}
1&\it when& x_0\not\in C_i\ \ \\
g\hbox{ \rm (mod  2)}&\it when& x_0\in C_i\ .
\end{array}\right.
$$
\end{pr}

{\ }
\vspace{10mm}  \\
{\small Christian Okonek: \\
Institut f\"ur Mathematik, Universit\"at Z\"urich,
Winterthurerstrasse 150, CH-8057 Z\"urich,\\
e-mail: okonek@math.unizh.ch
\\  \\
Andrei Teleman: \\
CMI,   Universit\'e de Provence,  39  Rue F. Joliot-Curie, F-13453
Marseille Cedex 13,   e-mail: teleman@cmi.univ-mrs.fr
}

  \end{document}